\documentclass[a4paper,leqno]{amsart}
\usepackage{graphicx}
\usepackage{amsmath}
\usepackage{amsfonts}

\newtheorem{theorem}{Theorem}[section]
\theoremstyle{plain}

\newtheorem{definition}{Definition}[section]

\numberwithin{equation}{section}
\numberwithin{equation}{section}

\begin{document}
\title[homogenization and imperfect contact]{Steady fluid flow in fissured
porous media with imperfect contact}
\author{Abdelhamid AINOUZ}
\address{Laboratory AMNEDP, Faculty of Mathematics, University of sciences
and technology, Po Box 32 El Alia, Algiers, Algeria.}
\email{aainouz@usthb.dz}
\thanks{The paper is a part of the project N$%
{{}^\circ}%
$ B00002278 of the MESRS algerian office. The author is indebted for their
financial support.}
\date{}

\begin{abstract}
In this paper, periodic homogenization of a steady fluid flow in fissured
porous solids with imperfect interfacial contact is performed via two-scale
asymptotic method.
\end{abstract}

\maketitle

\section{Introduction\label{s1}}

Micro-models for the fluid flow in fissured porous medium\ made of two
interacting porous systems with different permeabilities are very important
in reservoir petroleum, civil engineering, geophysics and many other areas
of engineering. In this paper, we shall deal with the homogenization of
fluid flow in fissured porous media. There are many papers devoted to the
subject, see for instance \cite{adh, adh1, cs} and the references therein.
Here, we shall be concerned with periodic porous media made of two
components with highly contrast of permeabilities and\ an imperfect contact
on the interface of these solids. One of the two porous structures is
associated with the fissures and the other one with the porous block or
matrix. More precisely, we shall consider a Newtonian fluid flow, in a
periodic fissured porous medium $\Omega =\Omega _{1}^{\varepsilon }\cup
\Omega _{2}^{\varepsilon }\cup \Gamma ^{\varepsilon }$, where $\Omega
_{1}^{\varepsilon }$ is the fissures region, $\Omega _{2}^{\varepsilon }$ is
the block and $\Gamma ^{\varepsilon }$ is the interface that separates these
two regions. The micro-model is based on mass conservation for the fluid in
each phase, combined with the Darcy's law . We assume that $\Omega
_{1}^{\varepsilon }\ $and $\Omega _{2}^{\varepsilon }$ are in imperfect
contact. That is, an exchange flow barrier formulation on the interface $%
\Gamma ^{\varepsilon }$ will be considered, see for instance, \cite{ain, ds,
ep}. The macro-model is derived by means of two-scale asymptotic technique (%
\textit{Cf}. \cite{blp, san}) and jusitified with the help of the two-scale
convergence method \cite{all}. It will be seen that the macro-model is, in
some sense, the limit of a family of periodic micro-models in which the size
of the periodicity approaches zero. It is shown that the overall behavior of
fluid flow in such micro-scale heterogeneous media, is a standard one.
Namely, the macroscopic equation is of elliptic type and it has the form $-$%
\textrm{div}$\left( A^{h}\nabla u\right) =F$. The only novelty here is that
the fluid flow presents an extra source surface density of the exterior
boundary. We mention that no double porosity effects occur at the macroscale
description.

The paper is organized as follows: Section 2 is devoted to the problem
setting of the micro-model. In Section \ref{s3}, we shall be concerned with
the derivation of the homogenized model via the formal procedure of
two-scale asymptotic ananlysis. Section \ref{s4} is devoted to the rigorous
derivation of the homogenized model obtained in the previous section by the
two scale convergence method. Finally, we end this paper with some comments
and conclusion.

\section{Setting of the Problem\label{s2} and the main result}

We consider $\Omega $\ a bounded, connected and smooth open subset of $%
\mathbb{R}^{N}\ $($N\geq 2$). We assume that $\Omega $ is made of a large
number of identical cells $\varepsilon Y$, where $\varepsilon >0$ is a
sufficiently small parameter ($\varepsilon \ll 1$) and $Y=]0,1[^{N}$\
denotes the generic cell of periodicity. Let $Y$\ be divided as follows: $%
Y:=Y_{1}\cup Y_{2}\cup \Gamma $ where $Y_{1},$\ $Y_{2}$\ are two sub-domains
of $Y$ and $\Gamma $ is the interface between $Y_{1}$ and\ $Y_{2}$. Namely $%
\Gamma =\overline{Y_{1}}\cap \overline{Y_{2}}$. We denote $\nu $ the unit
normal of $\Gamma $, outward to $Y_{1}$. Let $\chi _{1}$ and $\chi _{2}$\
denote respectively the characteristic function of $Y_{1}$ and $Y_{2}$,
extended by $Y$-periodicity to $\mathbb{R}^{N}$. For $\ i=1$ or $2$ we set%
\begin{equation*}
\Omega _{\varepsilon }^{i}:=\{x\in \Omega :\chi _{i}(\frac{x}{\varepsilon }%
)=1\}.
\end{equation*}%
Let $\Gamma ^{\varepsilon }:=\overline{\Omega _{1}^{\varepsilon }}\cap
\overline{\Omega _{2}^{\varepsilon }}$.\ We will assume that $\overline{%
\Omega _{2}^{\varepsilon }}\subset \Omega $, so that $\partial \Omega
_{2}^{\varepsilon }=\Gamma ^{\varepsilon }$ and $\partial \Omega
_{1}^{\varepsilon }=\partial \Omega \cup \Gamma ^{\varepsilon }$. Let $%
Z_{i}=\cup _{k\in \mathbb{Z}^{N}}\left( Y_{i}+k\right) $. As in \cite{all},
we shall assume that\ $Z_{1}$ is smooth and connected open subset of $%
\mathbb{R}^{N}$.

Let $A$ and $B$ denote respectively the permeability of the medium $Z_{1}$
and $Z_{2}$. We assume that $A$ (resp. $B$) is continuous on $\mathbb{R}^{N}$%
, $Y-$periodic and satisfies the ellipticity condition:\textbf{\ }%
\begin{equation*}
A\xi \cdot \xi \geq C|\xi |^{2},\text{(resp. }B\xi \cdot \xi \geq C|\xi |^{2}%
\text{)}\ \ \xi \in \mathbb{R}^{N}
\end{equation*}%
where, here and in what follows, $C$ denotes various positive constants
independent of $\varepsilon $. Let $f_{i}$ be a measurable function
representing the internal source density of the fluid flow in $\Omega
_{i}^{\varepsilon }$ and let $h$ be the \textit{non-rescaled} hydraulic
permeability of the thin layer $\Gamma ^{\varepsilon }$. We suppose that $%
f_{i}\in L^{2}\left( \Omega \right) $ $\ $and that $h$ is a continous
function on $\mathbb{R}^{N}$, $Y-$periodic and bounded from below:%
\begin{equation*}
h\left( y\right) \geq h_{0}>0\text{ }x\in \mathbb{R}^{N}.
\end{equation*}

To deal with periodic homogenization with microstructures, we shall denote
for $x\in \mathbb{R}^{N}$, $A^{\varepsilon }\left( x\right) =A\left( \frac{x%
}{\varepsilon }\right) $, $B^{\varepsilon }\left( x\right) =B\left( \frac{x}{%
\varepsilon }\right) $ and $h^{\varepsilon }\left( x\right) =\varepsilon
h\left( \frac{x}{\varepsilon }\right) $.

The micro-model that we shall study is given by the following set of
equations:
\begin{subequations}
\begin{align}
& -\mathrm{div}\left( A^{\varepsilon }\nabla u^{\varepsilon }\right) =f_{1}%
\text{ in }\Omega _{1}^{\varepsilon }\text{,}  \label{1b} \\
& -\varepsilon ^{2}\mathrm{div}\left( B^{\varepsilon }\nabla v^{\varepsilon
}\right) =f_{2}\text{ in }\Omega _{2}^{\varepsilon }\text{,}  \label{1a} \\
& A^{\varepsilon }\nabla u^{\varepsilon }\cdot \nu ^{\varepsilon
}=\varepsilon ^{2}B^{\varepsilon }\nabla v^{\varepsilon }\cdot \nu
^{\varepsilon }\text{ on }\Gamma ^{\varepsilon }\text{,}  \label{1c} \\
& A^{\varepsilon }\nabla u^{\varepsilon }\cdot \nu ^{\varepsilon
}=-\varepsilon h^{\varepsilon }\left( u^{\varepsilon }-v^{\varepsilon
}\right) \text{ on }\Gamma ^{\varepsilon }\text{,}  \label{1d} \\
& u^{\varepsilon }=0\text{ on }\partial \Omega \   \label{1e}
\end{align}%
where $\nu ^{\varepsilon }$ stands for the unit normal of $\Gamma
^{\varepsilon }$ outward to $\Omega _{1}^{\varepsilon }$. It is obtained by
the $Y-$periodicity extension of $\nu $. Let us mention that $\Omega
_{1}^{\varepsilon }$ represents the fissured space region with permeability $%
A$ and $\Omega _{2}^{\varepsilon }$ the block region with permeability $%
\varepsilon ^{2}B$. The quantitities $u^{\varepsilon }$ and $v^{\varepsilon
} $ are the fluid flow velocities in $\Omega _{1}^{\varepsilon }$ and $%
\Omega _{2}^{\varepsilon }$ respectively. Note tha we have chosen a
particular scaling of the permeability coefficients in (\ref{1a}). This
means that both terms $\int_{\Omega _{1}^{\varepsilon }}\left\vert \nabla
u^{\varepsilon }\right\vert ^{2}dx$ and $\varepsilon ^{2}\int_{\Omega
_{2}^{\varepsilon }}\left\vert \nabla v^{\varepsilon }\right\vert ^{2}dx$
have the same order of magnitude and thus leading to a balance in potential
energies. Equations (\ref{1b}) and (\ref{1a}) express the conservation of
mass of fluid with Darcy's law in $\Omega _{1}^{\varepsilon }$ and $\Omega
_{2}^{\varepsilon }$ respectively. The first equation describes the flow in
the fissured regions with large permeability and the second describes the
flow in the block system regions with low permeability. For more details, we
refer the reader to Arbogast, Douglas, and Hornung \cite{adh1} (see also
Allaire \cite{all}). Condition (\ref{1c}) expresses flux continuity across $%
\Gamma ^{\varepsilon } $, the interfacial flow thin layer with permeability
given by $h^{\varepsilon }$. Condition (\ref{1d}) models the imperfect
contact between the block and the fissures \cite{ds, ep}. The condition (\ref%
{1e}) is the standard homogeneous Dirichlet condition on the exterior
boundary of $\Omega $.

To set the mathematical framework of our Porblem (\ref{1b})-(\ref{1e}), we
need to introduce the following space $H^{\varepsilon }=\left( H^{1}\left(
\Omega _{1}^{\varepsilon }\right) \cap H_{0}^{1}\left( \Omega \right)
\right) \times H^{1}\left( \Omega _{2}^{\varepsilon }\right) $. The space $%
H^{\varepsilon }$ is equipped with the norm
\end{subequations}
\begin{equation*}
\left\Vert \left( \varphi ,\psi \right) \right\Vert _{H^{\varepsilon
}}^{2}=\left\Vert \nabla \varphi \right\Vert _{L^{2}\left( \Omega
_{1}^{\varepsilon }\right) }^{2}+\varepsilon ^{2}\left\Vert \nabla \psi
\right\Vert _{L^{2}\left( \Omega _{2}^{\varepsilon }\right)
}^{2}+\varepsilon \left\Vert \varphi -\psi \right\Vert _{L^{2}(\Gamma
^{\varepsilon })}^{2}\text{.}
\end{equation*}%
The weak formulation of (\ref{1b})-(\ref{1e}) can be read as follows\textbf{%
: }find $\ \left( u^{\varepsilon },v^{\varepsilon }\right) \in
H^{\varepsilon }$, such that for all $v=\left( \varphi ,\psi \right) \in
H^{\varepsilon }$, we have%
\begin{eqnarray}
&&\int_{\Omega _{1}^{\varepsilon }}A^{\varepsilon }\nabla v^{\varepsilon
}\nabla \varphi dx+\varepsilon ^{2}\int_{\Omega _{2}^{\varepsilon
}}B^{\varepsilon }\nabla v^{\varepsilon }\nabla \psi dx+  \notag \\
&&\int_{\Gamma ^{\varepsilon }}h^{\varepsilon }\left( v^{\varepsilon
}-v^{\varepsilon }\right) \left( \varphi -\psi \right) ds^{\varepsilon
}\left( x\right) =\int_{\Omega _{1}^{\varepsilon }}f_{1}\varphi
dx+\int_{\Omega _{2}^{\varepsilon }}f_{2}\psi dx\text{.}  \label{fv}
\end{eqnarray}%
where $dx$ and $ds^{\varepsilon }\left( x\right) $ denote respectively the
Lebesgue measure on $\mathbb{R}^{N}$ and the Hausdorff measure on $\Gamma
^{\varepsilon }$.

In view of the assumptions made on $A_{i}$, $f_{i}$ and $h$, we can easily
establish the following existence and uniqueness result whose proof is a
slight modification of that\ given by H. I. Ene and D. Polisevski \cite{ep}
and therefore will be omitted.

\begin{theorem}
For any sufficiently small $\varepsilon >0$, there exists a unique couple $%
\left( u^{\varepsilon },v^{\varepsilon }\right) \in H^{\varepsilon }$,
solution of the weak problem (\ref{fv}), such that
\begin{equation}
\Vert \left( u^{\varepsilon },v^{\varepsilon }\right) \Vert _{H^{\varepsilon
}}\leq C\text{.}  \label{30}
\end{equation}
\end{theorem}

Now, thanks to the a priori estimates (\ref{30}), one is led to study the
limiting behavior of\ the sequence $\left( u^{\varepsilon },v^{\varepsilon
}\right) $ as $\varepsilon $ approaches $0$. This is summarized in the main
result of the paper:

\begin{theorem}
\label{thp}Let $\left( u^{\varepsilon },v^{\varepsilon }\right) \in
H^{\varepsilon }$ be the solution of the weak system (\ref{fv}). Then, up to
a subsequence, there exists $u\in H_{0}^{1}\left( \Omega \right) $ and $%
v_{0}\in L^{2}\left( \Omega ;H_{\#}^{1}\left( Y_{2}\right) \right) $ such
that
\begin{eqnarray*}
u^{\varepsilon } &\rightharpoonup &u\text{ in }L^{2}\left( \Omega \right)
\text{ weakly,} \\
v^{\varepsilon } &\rightharpoonup &\int_{Y_{2}}v_{0}\left( y\right) dy\text{
in }L^{2}\left( \Omega \right) \text{ weakly,}
\end{eqnarray*}%
Let $w^{\varepsilon }=\chi _{1}\left( \frac{x}{\varepsilon }\right)
u^{\varepsilon }+\chi _{2}\left( \frac{x}{\varepsilon }\right)
v^{\varepsilon }$ denote the overall pressure. Then, $w^{\varepsilon }$
weakly converges to $U=u+G$ where $U$ is the unique solution to the
homogenized model:%
\begin{eqnarray*}
&&-\mathrm{div}\left( A^{\mathrm{h}}\nabla U\right) \text{%
$=$%
}F\text{ in }\Omega , \\
&&U\text{%
$=$%
}G\text{ on }\partial \Omega \text{.}
\end{eqnarray*}

Here, $A^{\mathrm{h}}$ (resp. $F,$ $G$) is given by (\ref{a}) (resp. (\ref%
{hs2}),(\ref{g})).
\end{theorem}

To prove this theorem, we shall first use the two-scale asymptotic procedure
to formally derive the homogenized model and apply then the two-scale
convergence technique to rigorously justify the homogenization. This is the
scope of the next two Sections.

\section{The formal homogenization procedure\label{s3}}

The purpose of this section is to formally construct the homogenized system
of (\ref{1b})-(\ref{1e}), via the two-scale asymptotic expansion method \cite%
{blp, san}. Assume the following ansatz:
\begin{eqnarray}
u^{\varepsilon }\left( x\right) &=&u_{0}\left( x,y\right) +\varepsilon
u_{1}\left( x,y\right) +\varepsilon ^{2}u_{2}\left( x,y\right) +\ldots ,
\label{ex1} \\
v^{\varepsilon }\left( x\right) &=&v_{0}\left( x,y\right) +\varepsilon
v_{1}\left( x,y\right) +\varepsilon ^{2}v_{2}\left( x,y\right) +\ldots \ \
\label{ex2}
\end{eqnarray}%
where $y=x/\varepsilon $, the unknowns $u_{0},\ v_{0},\ u_{1},v_{1},\
u_{2},v_{2}...$ are $Y$-periodic in the second variable $y$. Plugging the
above expansions (\ref{ex1})-(\ref{ex2}) into the set of equations (\ref{1b}%
--\ref{1e}) and identifying powers of $\varepsilon $ yields a hierarchy of
boundary value problems. At the first step, it is not difficult to observe
that Equation (\ref{1a}) at $\varepsilon ^{-2}$, and Equation (\ref{1c}) at $%
\varepsilon ^{-1}$ orders yield that $u_{0}$ is independent of $y\in Y_{1}$.
That is
\begin{equation*}
u_{0}\left( x,y\right) =u\left( x\right) .
\end{equation*}%
Next, Equation (\ref{1a}) at $\varepsilon ^{-1}$ and Equation (\ref{1c}) at $%
\varepsilon ^{0}$ orders show that the corrector term $u_{1}$ may be written
as
\begin{equation}
u_{1}\left( x,y\right) =\sum_{j=1}^{N}\frac{\partial u}{\partial x_{j}}%
\left( x\right) \omega _{j}\left( y\right) +\tilde{u}\left( x\right)
\label{3m}
\end{equation}%
where, for $1\leq j\leq N$, $\omega _{j}\in H_{\#}^{1}\left( Y_{1}\right) /%
\mathbb{R}$ is the unique solution to the following cell problem:
\begin{equation*}
\int_{Y_{1}}\left( A\nabla _{y}\omega _{j},\nabla _{y}\zeta \right)
dy=\int_{Y_{1}}\left( -Ae_{j},\nabla _{y}\zeta \right) dy,\ \zeta \in
H_{\#}^{1}\left( Y_{1}\right) \text{,}
\end{equation*}%
where $\left( e_{j}\right) $ is the canonical basis of $\mathbb{R}^{N}$. In (%
\ref{3m}), $\tilde{u}\left( x\right) $ is any additive constant. At the
final step, Equation (\ref{1b})-(\ref{1a}) at $\varepsilon ^{0}$ and
Equations (\ref{1c})-(\ref{1d})) at $\varepsilon ^{1}$ orders yield the
system:
\begin{equation}
-\mathrm{div}_{y}\left( A\nabla _{y}u_{2}\right) =f_{1}+\mathrm{div}%
_{y}\left( A\nabla _{x}u_{1}\right) +\mathrm{div}_{x}\left( A\left( \nabla
_{y}u_{1}+\nabla _{x}u\right) \right) \ \text{in }\Omega \times Y_{1},
\label{4o}
\end{equation}%
\begin{equation}
-\mathrm{div}_{y}\left( B\nabla _{y}v_{0}\right) =f_{2}\text{ in }\Omega
\times Y_{2};  \label{4n}
\end{equation}%
\begin{equation}
A\nabla _{y}u_{2}\cdot \nu =-A\nabla _{x}u_{1}\cdot \nu +B\nabla
_{y}v_{0}\cdot \nu \text{ on }\Omega \times \Gamma ,  \label{4p}
\end{equation}%
\begin{equation}
A\nabla _{y}u_{2}\cdot \nu =-A\nabla _{x}u_{1}\cdot \nu -h\left(
u-v_{0}\right) \text{ on }\Omega \times \Gamma ,  \label{4q}
\end{equation}%
\begin{equation}
y\mapsto v_{0}\left( x,y\right) ,\ u_{2}\left( x,y\right) \text{ }Y-\text{%
periodic.}  \label{4r}
\end{equation}

The weak formulation of Equations (\ref{4o}), (\ref{4q})-(\ref{4r}) is
\begin{equation*}
\int_{Y_{1}}\left( A\nabla _{y}u_{2},\nabla _{y}\zeta \right)
dy=\left\langle F,\zeta \right\rangle ,\ \ \zeta \in H_{\#}^{1}\left(
Y_{1}\right) \text{,}
\end{equation*}%
\ where
\begin{eqnarray*}
\left\langle F,\zeta \right\rangle &=&\int_{Y_{1}}\mathrm{div}_{x}\left(
A\left( \nabla _{y}u_{1}+\nabla u\right) \right) \zeta \\
&&-\int_{Y_{1}}A\nabla _{x}u_{1}\nabla _{y}\zeta +\int_{\Gamma }h\left(
u-v_{0}\right) \zeta ds\left( y\right) +\int_{Y_{1}}f_{1}\zeta
\end{eqnarray*}%
where $ds\left( y\right) $ denotes the Hausdorff measure on $\Gamma $.
Again, using the Divergence Theorem ( as in \cite{san} ), a necessary and
sufficient condition for the existence and uniqueness of $u_{2}\in
H_{\#}^{1}\left( Y_{1}\right) /\mathbb{R}$ is that $F$ satisfies the
compatibility condition $\left\langle F,1\right\rangle =0$. This reads in
our case as follows:
\begin{equation}
\int_{Y_{1}}\left( \mathrm{div}_{x}\left( A\left( \nabla _{y}u_{1}+\nabla
u\right) \right) \right) dy+\int_{\Gamma }h\left( u-v_{0}\right) ds\left(
y\right) =\int_{Y_{1}}f_{1}dy.  \label{5a}
\end{equation}%
Using (\ref{3m}), Equation (\ref{5a}) becomes
\begin{equation}
-\mathrm{div}\left( A^{\mathrm{h}}\nabla u\right) +\int_{\Gamma }h\left(
u-v_{0}\right) ds\left( y\right) =\int_{Y_{1}}f_{1}dy  \label{hs1}
\end{equation}%
where $A^{\mathrm{h}}=\left( a_{ij}^{\mathrm{h}}\right) _{1\leq i,j\leq N}$
is given by
\begin{equation}
a_{ij}^{\mathrm{h}}=\int_{Y_{1}}A\left( \nabla _{y}\omega _{i}+e_{i}\right)
\cdot \left( \nabla _{y}\omega _{j}+e_{j}\right) dy\text{.}  \label{a}
\end{equation}%
On the other hand, it is easy to see that $v_{0}$ can be written as
\begin{equation}
v_{0}\left( x,y\right) -u\left( x\right) =\alpha \left( y\right) f_{2}\left(
x\right) ,\ \ \left( x,y\right) \in \Omega \times Y_{2}  \label{4m}
\end{equation}%
where $\alpha \in H_{\#}^{1}\left( Y_{2}\right) $ is the unique solution of
the following problem
\begin{equation}
\left\{
\begin{array}{l}
-\mathrm{div}_{y}\left( B\nabla _{y}\alpha \right) =1\text{ in }Y_{2}, \\
\\
B\nabla _{y}\alpha \cdot \nu -h\alpha =0\text{ on }\Gamma .%
\end{array}%
\right.  \label{alpha}
\end{equation}%
Therefore, Equation (\ref{hs1}) yields
\begin{equation}
-\mathrm{div}\left( A^{\mathrm{h}}\nabla u\right) =|Y_{1}|f_{1}+\hat{\alpha}%
f_{2}:=F^{\ast }  \label{hs2}
\end{equation}%
where $\hat{\alpha}=\int_{\Gamma }\alpha ds\left( y\right) $. The boundary
condition for $u$ is obtained from (\ref{1e}) at $\varepsilon ^{0}$ order
and it reads
\begin{equation}
u=0\text{ on }\partial \Omega \text{.}  \label{b1}
\end{equation}

Let us observe that the overall pressure $w^{\varepsilon }=\chi _{1}\left(
\frac{x}{\varepsilon }\right) u^{\varepsilon }+\chi _{2}\left( \frac{x}{%
\varepsilon }\right) v^{\varepsilon }$ two scale converges to $u+\chi
_{2}\alpha f_{2}$ (see Definition \ref{def1} below). Consequently, $%
w^{\varepsilon }$ weakly converges to $U=u+G$ where
\begin{equation}
G=\left( \int_{Y_{2}}\alpha \right) f_{2}.  \label{g}
\end{equation}
Equation (\ref{hs2}) is the so-called macroscopic equation for $u$.
Moreover, if $G$ is sufficiently smooth, then the homogenized model for the
weak limit $U$ is as follows:%
\begin{eqnarray}
&&-\mathrm{div}\left( A^{\mathrm{h}}\nabla U\right) =F\text{ in }\Omega ,
\label{ndhp1} \\
&&U=G\text{ on }\partial \Omega  \label{ndhp2}
\end{eqnarray}%
where $F=F^{\ast }+\mathrm{div}\left( A^{\mathrm{h}}\nabla G\right) $. Note
that in the homogenized problem (\ref{ndhp1})-(\ref{ndhp2}), a non
homogeneous Dirichlet is derived, on the contrary to the micro-model, where
homogeneous boundary condition is prescribed, see condition (\ref{1e}). This
extra source surface density essentially arises from the fact that

\begin{enumerate}
\item Blocks have low permeability;

\item non null source volumetric density on the blocks;

\item and the contribution of the fluid flow in the vicinity of the thin
layer where the contact is assumed imperfect.
\end{enumerate}

\section{Two-scale convergence approach \label{s4}}

In this section, we will derive the homogenized system of (\ref{1b})-(\ref%
{1e}) by the two scale convergence method, see \cite{all}. First, we define $%
C_{\#}(Y)$ to be the space of all continuous functions on $\mathbb{R}^{3}$
which are $Y$-periodic. Let the space $L_{\#}^{2}\left( Y\right) $ (resp. $%
L_{\#}^{2}\left( Y_{i}\right) $, $i=1,2$) to be all functions belonging to $%
L_{\mathrm{loc}}^{2}\left( \mathbb{R}^{3}\right) $ (resp. $L_{\mathrm{loc}%
}^{2}\left( Z_{i}\right) $) which are $Y$-periodic, and $H_{\#}^{1}\left(
Y\right) $ (resp. $H_{\#}^{1}\left( Y_{i}\right) $) to be the space of those
functions together with their derivatives belonging to $L_{\#}^{2}\left(
Y\right) $ (resp. $L_{\#}^{2}\left( Z_{i}\right) $). Now, we recall the
definition and main results concerning the method of two-scale convergence.
For full details, we refer the reader to \cite{all, adh}.

\begin{definition}
\label{def1}A sequence $\left( v^{\varepsilon }\right) \ $in $L^{2}\left(
\Omega \right) $ two-scale converges to $v\in L^{2}\left( \Omega \times
Y\right) $ (we write $v^{\varepsilon }\overset{2-s}{\rightharpoonup }v$) if,
for any admissible test function $\varphi \in L^{2}\left( \Omega ;\mathcal{C}%
_{\#}(Y)\right) $,
\begin{equation*}
\lim_{\varepsilon \rightarrow 0}\int_{\Omega }v^{\varepsilon }\left(
x\right) \varphi \left( x,\frac{x}{\varepsilon }\right) dx=\int_{\Omega
\times Y}v\left( x,y\right) \varphi \left( x,y\right) dxdy\text{.}
\end{equation*}
\end{definition}

\begin{theorem}
\label{t1}Let $(v^{\varepsilon })$ be a sequence of functions in $%
L^{2}(\Omega )$ which is uniformly bounded. Then, there exist $v\in
L^{2}(\Omega \times Y)$ and a subsequence of $(v^{\varepsilon })$ which
two-scale converges\ to $v$.
\end{theorem}

\begin{theorem}
\label{t2}Let $(v^{\varepsilon })$ be a uniformly bounded sequence in $%
H^{1}(\Omega )$ (resp. $H_{0}^{1}(\Omega )$). Then there exist $v\in
H^{1}\left( \Omega \right) $ (resp. $H_{0}^{1}(\Omega )$) and $\hat{v}\in
L^{2}(\Omega ;H_{\#}^{1}(Y)/\mathbb{R})$ such that, up to a subsequence,%
\begin{equation*}
v^{\varepsilon }\overset{2-s}{\rightharpoonup }v;\ \ \ \ \ \nabla
v^{\varepsilon }\overset{2-s}{\rightharpoonup }\nabla v+\nabla _{y}\hat{v}.
\end{equation*}
\end{theorem}

\begin{theorem}
\label{t3}Let $(v^{\varepsilon })$ be a sequence of functions in $%
H^{1}(\Omega )$ such that
\begin{equation*}
\left\Vert v^{\varepsilon }\right\Vert _{L^{2}\left( \Omega \right)
}+\varepsilon \left\Vert \nabla v^{\varepsilon }\right\Vert _{L^{2}\left(
\Omega \right) ^{3}}\leq C\text{.}
\end{equation*}%
Then, there exist $v\in L^{2}\left( \Omega ;H_{\#}^{1}(Y)\right) $ and a
subsequence of $\left( v^{\varepsilon }\right) $, still denoted by $\left(
v^{\varepsilon }\right) $ such that
\begin{equation*}
v^{\varepsilon }\overset{2-s}{\rightharpoonup }v,\ \ \ \ \ \varepsilon
\nabla v^{\varepsilon }\overset{2-s}{\rightharpoonup }\nabla _{y}v
\end{equation*}%
and for every $\varphi \in \mathcal{D}\left( \Omega ;\mathcal{C}%
_{\#}(Y)\right) $, we have:%
\begin{equation*}
\lim_{\varepsilon \rightarrow 0}\int_{\Gamma ^{\varepsilon }}\varepsilon
v^{\varepsilon }\left( x\right) \varphi \left( x,\frac{x}{\varepsilon }%
\right) ds^{\varepsilon }\left( x\right) =\int_{\Omega \times \Gamma
}v\left( x,y\right) \varphi \left( x,y\right) dxds\left( y\right) .
\end{equation*}%
Here and in the sequel $ds\left( y\right) $ denotes the Hausdorff measure on
$\Gamma $.
\end{theorem}

As a direct application of the theorems listed above (Thms \ref{t1}-\ref{t3}%
) and the a priori estimates (\ref{30}), we give without proof the following
two-scale convergence result concerning the solutions $\left( u^{\varepsilon
},v^{\varepsilon }\right) $ of the Problem (\ref{fv}).

\begin{theorem}
\label{t4}There exists a subsequence of $\left( u^{\varepsilon
},v^{\varepsilon }\right) $, solution of (\ref{fv}), still denoted $\left(
u^{\varepsilon },v^{\varepsilon }\right) $, and there exist unique $u\in
H_{0}^{1}\left( \Omega \right) $,\ \ $\ u_{1}\in L^{2}(\Omega
;H_{\#}^{1}\left( Y\right) /\mathbb{R})$ and $v_{0}\in L^{2}\left( \Omega
;H_{\#}^{1}\left( Y_{2}\right) \right) $ such that
\begin{equation}
\chi _{1}^{\varepsilon }u^{\varepsilon }\overset{2-s}{\rightharpoonup }\chi
_{1}u,\ \chi _{2}^{\varepsilon }v^{\varepsilon }\overset{2-s}{%
\rightharpoonup }\chi _{2}v_{0}  \label{321}
\end{equation}%
\ and
\begin{equation}
\chi _{1}^{\varepsilon }\nabla u^{\varepsilon }\overset{2-s}{\rightharpoonup
}\chi _{1}\left( \nabla u+\nabla _{y}u_{1}\right) ,\ \varepsilon \chi
_{2}^{\varepsilon }\nabla v^{\varepsilon }\overset{2-s}{\rightharpoonup }%
\chi _{2}\nabla _{y}v_{0}.  \label{35}
\end{equation}%
Moreover, the following convergence holds: \
\begin{equation}
\lim_{\varepsilon \rightarrow 0}\int_{\Gamma ^{\varepsilon }}\varepsilon
\left( u^{\varepsilon }-v^{\varepsilon }\right) \psi ^{\varepsilon
}ds^{\varepsilon }\left( x\right) =\int_{\Omega \times \Gamma }\left(
u-v_{0}\right) \psi dxds\left( y\right) ,  \label{36}
\end{equation}%
for any $\psi \in \mathcal{D}\left( \Omega ;\mathcal{C}_{\#}\left( Y\right)
\right) $ with $\psi ^{\varepsilon }\left( x\right) =\psi \left(
x,x/\varepsilon \right) $.
\end{theorem}

To determine the limiting equations of the system (\ref{fv}), we begin by
choosing the adequate admissible test functions. Let $\varphi ^{\varepsilon
}(x)=\varphi \left( x\right) +\varepsilon \varphi _{1}\left( x,\dfrac{x}{%
\varepsilon }\right) $ and $\psi ^{\varepsilon }(x)=\psi \left( x,\dfrac{x}{%
\varepsilon }\right) $ where $\varphi \in \mathcal{D}\left( \Omega \right) $
and $\varphi _{1}$, $\psi \in \mathcal{D}\left( \Omega ;\mathcal{C}%
_{\#}^{\infty }\left( Y\right) \right) $. Taking $\varphi =\varphi
^{\varepsilon }$ and $\psi =\psi ^{\varepsilon }$ in (\ref{fv}), we obtain
\begin{eqnarray}
&&\int_{\Omega _{1}^{\varepsilon }}A\left( \frac{x}{\varepsilon }\right)
\nabla u^{\varepsilon }\left( \nabla \varphi +\nabla _{y}\varphi _{1}\right)
dx+\int_{Q_{2}^{\varepsilon }}\varepsilon B\left( \frac{x}{\varepsilon }%
\right) \nabla v^{\varepsilon }\nabla _{y}\psi dx+  \notag \\
&&\varepsilon \int_{\Gamma ^{\varepsilon }}h\left( \frac{x}{\varepsilon }%
\right) \left( u^{\varepsilon }-v^{\varepsilon }\right) \left( \varphi -\psi
\right) ds^{\varepsilon }\left( x\right) +\varepsilon R^{\varepsilon
}=\int_{\Omega _{1}^{\varepsilon }}f_{1}\varphi dx+\int_{\Omega
_{2}^{\varepsilon }}f_{2}\psi dx  \label{46}
\end{eqnarray}%
where
\begin{eqnarray*}
&&R^{\varepsilon }=\int_{\Omega _{1}^{\varepsilon }}A\left( \frac{x}{%
\varepsilon }\right) \nabla u^{\varepsilon }\nabla _{x}\varphi _{1}\left( x,%
\frac{x}{\varepsilon }\right) dx+\varepsilon \int_{\Omega _{2}^{\varepsilon
}}B\left( \frac{x}{\varepsilon }\right) \nabla v^{\varepsilon }\nabla
_{x}\psi \left( x,\frac{x}{\varepsilon }\right) dx+ \\
&&\varepsilon \int_{\Gamma ^{\varepsilon }}h\left( \frac{x}{\varepsilon }%
\right) \left( u^{\varepsilon }-v^{\varepsilon }\right) \varphi _{1}\left( x,%
\frac{x}{\varepsilon }\right) ds^{\varepsilon }\left( x\right) \text{.}
\end{eqnarray*}%
In view of (\ref{35}), one can deduce that
\begin{eqnarray*}
&&\lim_{\varepsilon \rightarrow 0}\int_{\Omega _{1}^{\varepsilon }}A\left(
\frac{x}{\varepsilon }\right) \nabla v^{\varepsilon }\left( \nabla \varphi
\left( x\right) +\nabla _{y}\varphi _{1}\left( x,\frac{x}{\varepsilon }%
\right) \right) dx= \\
&&\int_{\Omega \times Y}\chi _{1}\left( y\right) A\left( y\right) \left(
\nabla u+\nabla _{y}u_{1}\right) \left( \nabla \varphi \left( x\right)
+\nabla _{y}\varphi _{1}\left( x,y\right) \right) dxdy
\end{eqnarray*}%
and
\begin{eqnarray*}
&&\lim_{\varepsilon \rightarrow 0}\int_{\Omega _{2}^{\varepsilon
}}\varepsilon B\left( \frac{x}{\varepsilon }\right) \nabla v^{\varepsilon
}\nabla _{y}\psi \left( x,\frac{x}{\varepsilon }\right) dx= \\
&&\int_{\Omega \times Y}\chi _{2}\left( y\right) B\left( y\right) \nabla
v_{0}\nabla _{y}\psi \left( x,y\right) dxdy.
\end{eqnarray*}%
By virtue of (\ref{36}), we find that
\begin{eqnarray*}
&&\lim_{\varepsilon \rightarrow 0}\varepsilon \int_{\Gamma ^{\varepsilon
}}h\left( \frac{x}{\varepsilon }\right) \left( u^{\varepsilon
}-v^{\varepsilon }\right) \left( \varphi \left( x\right) -\psi \left( x,%
\frac{x}{\varepsilon }\right) \right) ds^{\varepsilon }\left( x\right) = \\
&&\int_{\Omega \times \Gamma }h\left( y\right) \left( u-v_{0}\right) \left(
\varphi \left( x\right) -\psi \left( x,y\right) \right) dxds\left( y\right) .
\end{eqnarray*}%
We observe that $R^{\varepsilon }=O\left( 1\right) $ and, by collecting all
the preceeding limits, we get the following limiting equation of (\ref{fv}):%
\begin{eqnarray}
&&\int_{\Omega \times Y_{1}}A\left( y\right) \left( \nabla u+\nabla
_{y}u_{1}\right) \left( \nabla \varphi +\nabla _{y}\varphi _{1}\right)
dxdy+\int_{\Omega \times Y_{2}}B\left( y\right) \nabla _{y}v_{0}\nabla
_{y}\psi dxdy+  \notag \\
&&\int_{\Omega \times \Gamma }h\left( y\right) \left( u-v_{0}\right) \left(
\varphi -\psi \right) dxds\left( y\right)  \notag \\
&=&\int_{\Omega \times Y_{1}}f_{1}\varphi dx+\int_{\Omega \times
Y_{2}}f_{2}\psi dx\text{.}  \label{50}
\end{eqnarray}%
By density argument, the equation (\ref{50}) still holds true for any $%
\left( \varphi ,\varphi _{1},\varphi _{2}\right) \in H_{0}^{1}\left( \Omega
\right) \times L^{2}\left( \Omega ;H_{\#}^{1}\left( Y_{1}\right) /\mathbb{R}%
\right) \times L^{2}\left( \Omega ;H_{\#}^{1}\left( Y_{2}\right) \right) $.
We can summarize the preceding by observing that these equations are a weak
formulation associated to the two-scale homogenized system (\ref{474})-(\ref%
{482}) below:
\begin{eqnarray}
&&-\mathrm{div}_{y}\left( A\left( \nabla u+\nabla _{y}u_{1}\right) \right)
\text{%
$=$%
\ }0\text{ a.e. in }\Omega \times Y_{1},  \label{474} \\
&&-\mathrm{div}_{y}\left( B\nabla _{y}v_{0}\right) \text{%
$=$%
\ }f_{2}\text{ a.e.in }\Omega \times Y_{2},  \label{475} \\
&&-\mathrm{div}\left( \int_{Y_{1}}A\left( \nabla u+\nabla _{y}u_{1}\right)
dy\right) +  \notag \\
&&\int_{\Gamma }h\left( y\right) \left[ u-v_{0}\right] ds\left( y\right)
\text{%
$=$%
\ }f_{1}\text{ a.e. in }\Omega ,  \label{476}
\end{eqnarray}%
with the transmission and boundary conditions:%
\begin{eqnarray}
&&\left( A\left( \nabla u+\nabla _{y}u_{1}\right) \right) \cdot \nu \text{%
$=$%
\ }0\text{ a.e. on }\Omega \times \Gamma ,  \label{478} \\
&&B\nabla _{y}v_{0}\cdot v\text{%
$=$%
\ }-h\left( y\right) \left[ u-v_{0}\right] \text{ a.e. on }\Omega \times
\Gamma ,  \label{479} \\
&&y\longmapsto \ u_{1},v_{0}\text{ }Y\text{-periodic,}  \label{481} \\
&&u=0\text{ on }\partial \Omega \text{.}  \label{482}
\end{eqnarray}

Let us first note that Equations (\ref{474}) and (\ref{478}) lead to the
relation (\ref{3m}).\ Similarly, Equations (\ref{475}), (\ref{479}) and (\ref%
{481}) yield (\ref{4m}). Equation (\ref{476}) is the same as (\ref{5a}).
With these remarks in mind, one can easily recover the homogenization
procedure \ of the preceeding section and thus we have rigorously justified
the formal two-scale asymptotic expansion method.

\section{Conclusion}

We have used the homogenization theory to derive a macro-model for fluid
flow in fissured media with microstructures, in which blocks present very
low permeabilities. Moreover, the contact between the block region and the
fissures region is of imperfect type. The micro-model that we considered in
this paper is with a homogeneous Dirichlet boundary condition prescribed on
the exterior boundary. We have shown that the overall behavior of fluid flow
in such heterogeneous media with low permeability at the micro-scale is a
classical problem except that in the vicinity of the surface, there is an
additional source density arising from the source density of the blocks at
the micro-scale and the fact that the contact at micro-scale is imperfect.

\end{document}